\documentclass[9pt, a4paper]{amsart}

\usepackage{xcolor}
\definecolor{webgreen}{rgb}{0,.5,0}
\definecolor{webbrown}{rgb}{.6,0,0}
\definecolor{RoyalBlue}{cmyk}{1, 0.50, 0, 0}
\usepackage[colorlinks=true, breaklinks=true, urlcolor=webbrown, linkcolor=RoyalBlue, citecolor=webgreen,backref=page]{hyperref}

\usepackage{epsfig, graphicx, subfigure}
\usepackage{verbatim, setspace}
\usepackage{amsmath, amssymb}

\usepackage[T1]{fontenc}
\usepackage{textcomp}
\usepackage{fourier}

\newcommand{\D}		{\mathbb{D}}

\newcommand{\C}		{\mathbb{C}}
\newcommand{\N}		{\mathbb{N}}

\renewcommand{\det}{\mathrm{det}}

\newcommand{\qasq}{\quad \text{as} \quad}
\newcommand{\qandq}{\quad \text{and} \quad}

\newcommand{\rhy}   {\textnormal{RHP}-${\boldsymbol Y}$}
\newcommand{\rhx}   {\textnormal{RHP}-\( \boldsymbol X \)}
\newcommand{\rhm}   {\textnormal{RHP}-${\boldsymbol M}$}
\newcommand{\rhz}   {\textnormal{RHP}-${\boldsymbol Z}$}

\newcommand{\dd}{\mathrm{d}}
\newcommand{\ic}{\mathrm{i}}

\usepackage{varwidth}
\usepackage{varwidth}
\usepackage[most]{tcolorbox}

\tcbuselibrary{skins}
\tcbuselibrary{theorems}

\newenvironment{myitemize}{\begin{itemize}}{\end{itemize}}
\tcolorboxenvironment{myitemize}{blanker, breakable, before skip=6pt, after skip=6pt, borderline west={1mm}{0pt}{RoyalBlue}}

\newtcbtheorem{assumption}{Assumption}{enhanced,
before skip=2mm,after skip=2mm, colback=black!5,colframe=webgreen,boxrule=0.2mm,
attach boxed title to top left={xshift=1cm,yshift*=1mm-\tcboxedtitleheight}, varwidth boxed title*=-3cm,
boxed title style={frame code={
            \path[fill=tcbcolback!30!black]
              ([yshift=-1mm,xshift=-1mm]frame.north west)
                arc[start angle=0,end angle=180,radius=1mm]
              ([yshift=-1mm,xshift=1mm]frame.north east)
                arc[start angle=180,end angle=0,radius=1mm];
            \path[left color=tcbcolback!60!black,right color=tcbcolback!60!black,
              middle color=tcbcolback!80!black]
              ([xshift=-2mm]frame.north west) -- ([xshift=2mm]frame.north east)
              [rounded corners=1mm]-- ([xshift=1mm,yshift=-1mm]frame.north east)
              -- (frame.south east) -- (frame.south west)
              -- ([xshift=-1mm,yshift=-1mm]frame.north west)
              [sharp corners]-- cycle;
            },interior engine=empty,
          },
          fonttitle=\bfseries,
          title={#2},#1}{as}
          
 \newtcbtheorem{proposition}{Proposition}{enhanced,
before skip=2mm,after skip=2mm, colback=black!5,colframe=RoyalBlue,boxrule=0.2mm,
attach boxed title to top left={xshift=1cm,yshift*=1mm-\tcboxedtitleheight}, varwidth boxed title*=-3cm,
boxed title style={frame code={
            \path[fill=tcbcolback!30!black]
              ([yshift=-1mm,xshift=-1mm]frame.north west)
                arc[start angle=0,end angle=180,radius=1mm]
              ([yshift=-1mm,xshift=1mm]frame.north east)
                arc[start angle=180,end angle=0,radius=1mm];
            \path[left color=tcbcolback!60!black,right color=tcbcolback!60!black,
              middle color=tcbcolback!80!black]
              ([xshift=-2mm]frame.north west) -- ([xshift=2mm]frame.north east)
              [rounded corners=1mm]-- ([xshift=1mm,yshift=-1mm]frame.north east)
              -- (frame.south east) -- (frame.south west)
              -- ([xshift=-1mm,yshift=-1mm]frame.north west)
              [sharp corners]-- cycle;
            },interior engine=empty,
          },
          fonttitle=\bfseries,
          title={#2},#1}{prop}         

 \newtcbtheorem{theorem}{Theorem}{enhanced,
before skip=2mm,after skip=2mm, colback=black!5,colframe=webbrown,boxrule=0.2mm,
attach boxed title to top left={xshift=1cm,yshift*=1mm-\tcboxedtitleheight}, varwidth boxed title*=-3cm,
boxed title style={frame code={
            \path[fill=tcbcolback!30!black]
              ([yshift=-1mm,xshift=-1mm]frame.north west)
                arc[start angle=0,end angle=180,radius=1mm]
              ([yshift=-1mm,xshift=1mm]frame.north east)
                arc[start angle=180,end angle=0,radius=1mm];
            \path[left color=tcbcolback!60!black,right color=tcbcolback!60!black,
              middle color=tcbcolback!80!black]
              ([xshift=-2mm]frame.north west) -- ([xshift=2mm]frame.north east)
              [rounded corners=1mm]-- ([xshift=1mm,yshift=-1mm]frame.north east)
              -- (frame.south east) -- (frame.south west)
              -- ([xshift=-1mm,yshift=-1mm]frame.north west)
              [sharp corners]-- cycle;
            },interior engine=empty,
          },
          fonttitle=\bfseries,
          title={#2},#1}{thm}

\begin{document}

\title[Symmetric contours that separate the plane]{On multipoint Pad\'e approximants whose poles accumulate on contours that separate the plane}

\author{Maxim L. Yattselev}

\address{Department of Mathematical Sciences \\ Indiana University-Purdue University Indianapolis \\ 402~North Blackford Street, Indianapolis, IN 46202}

\address{Keldysh Institute of Applied Mathematics, Russian Academy of Science, Miusskaya Pl. 4, Moscow, 125047 Russian Federation}

\email{\href{mailto:maxyatts@iupui.edu}{maxyatts@iupui.edu}}

\thanks{The research was supported in part by a grant from the Simons Foundation, CGM-706591.}



\begin{abstract}
In this note we consider asymptotics of the multipoint Pad\'e approximants to Cauchy integrals of analytic non-vanishing densities defined on a Jordan arc connecting \( -1 \) and~\( 1 \). We allow for the situation where the (symmetric) contour attracting the poles of the approximants does separate the plane.
\end{abstract}

\maketitle

\section{Introduction and Main Results}

Let \( L \) be a smooth Jordan arc joining \( -1 \) and \( 1 \) (and oriented in this way), and
\[
w_L(z) := \sqrt{z^2-1}, \quad w_L(z) = z + \mathcal O(1), \qasq z\to\infty,
\]
be the branch holomorphic in \( \C\setminus L \). Let \( \rho(z) \) be a function analytic in a ``sufficiently large'', see Assumption~\ref{as:2} further below, neighborhood of 
\( L \) and
\begin{equation}
\label{f}
\widehat\rho_L(z) := \frac1{2\pi\ic}\int_L\frac{\rho(s)}{s-z}\frac{\dd s}{w_{L+}(s)}, \quad z\not\in L,
\end{equation}
where \( w_{L+}(s) \) is the trace of \( w_L(z) \) on the left-hand side of \( L \). Further, let \( \{E_i\}_{i=1}^\infty \) be an \emph{interpolation scheme}, i.e., a sequence of multi-sets in \( D_L :=\overline\C\setminus L \) such that each \( E_i \) consists of \( i \) not necessarily distinct nor finite points from \( D_L \). The \emph{multipoint Pad\'e approximant to \( \widehat\rho_L(z) \) of type \( (m,n) \) associated with \( E_{m+n} \)} is a rational function \( [m/n;E_{m+n}]_{\widehat\rho_L}(z) = (p_{m,n}/q_{m,n})(z) \) such that \( \deg(p_{m,n}) \leq m \), \( \deg(q_{m,n}) \leq n \), \( q_{m,n}(z)\not\equiv 0 \), and
\begin{equation}
\label{Rn}
\frac{q_{m,n}(z)\widehat\rho_L(z)-p_{m,n}(z)}{v_{m,n}(z)} = \mathcal O\left(z^{-\min\{m,n\}-1}\right) \qasq z\to\infty
\end{equation}
as well as is analytic in \( D_L \), where \( v_{m,n}(z) \) is the monic polynomial vanishing at the finite elements of \( E_{m+n} \) according to their multiplicity\footnote{This definition yields one additional interpolation condition at infinity.}, i.e.,
\[
v_{m,n}(z) := \prod_{e\in E_{m+n},|e|<\infty} (z-e).
\]
It is well-known that linear system \eqref{Rn} is always solvable and that the solution corresponding to monic \( q_{m,n}(z) \) of minimal degree is always unique. In what follows, we understand that \( q_{m,n}(z) \) and \( p_{m,n}(z) \) come from this unique solution. We shall call an approximant \emph{diagonal} if \( m=n \), in which case we simply write \( p_n(z) \) \( q_n(z) \), and \( v_n(z) \).

\begin{figure}[ht!]
\centering
\subfigure[]{\includegraphics[scale=0.4]{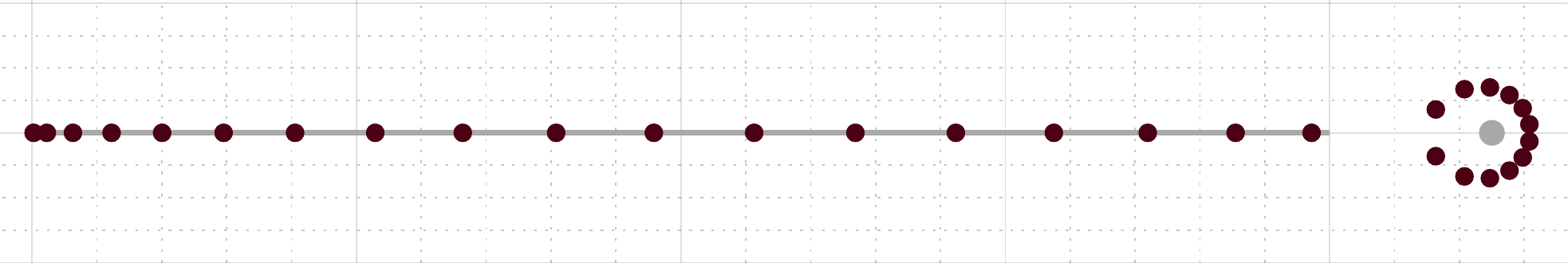}}
\quad
\subfigure[]{\includegraphics[scale=0.4]{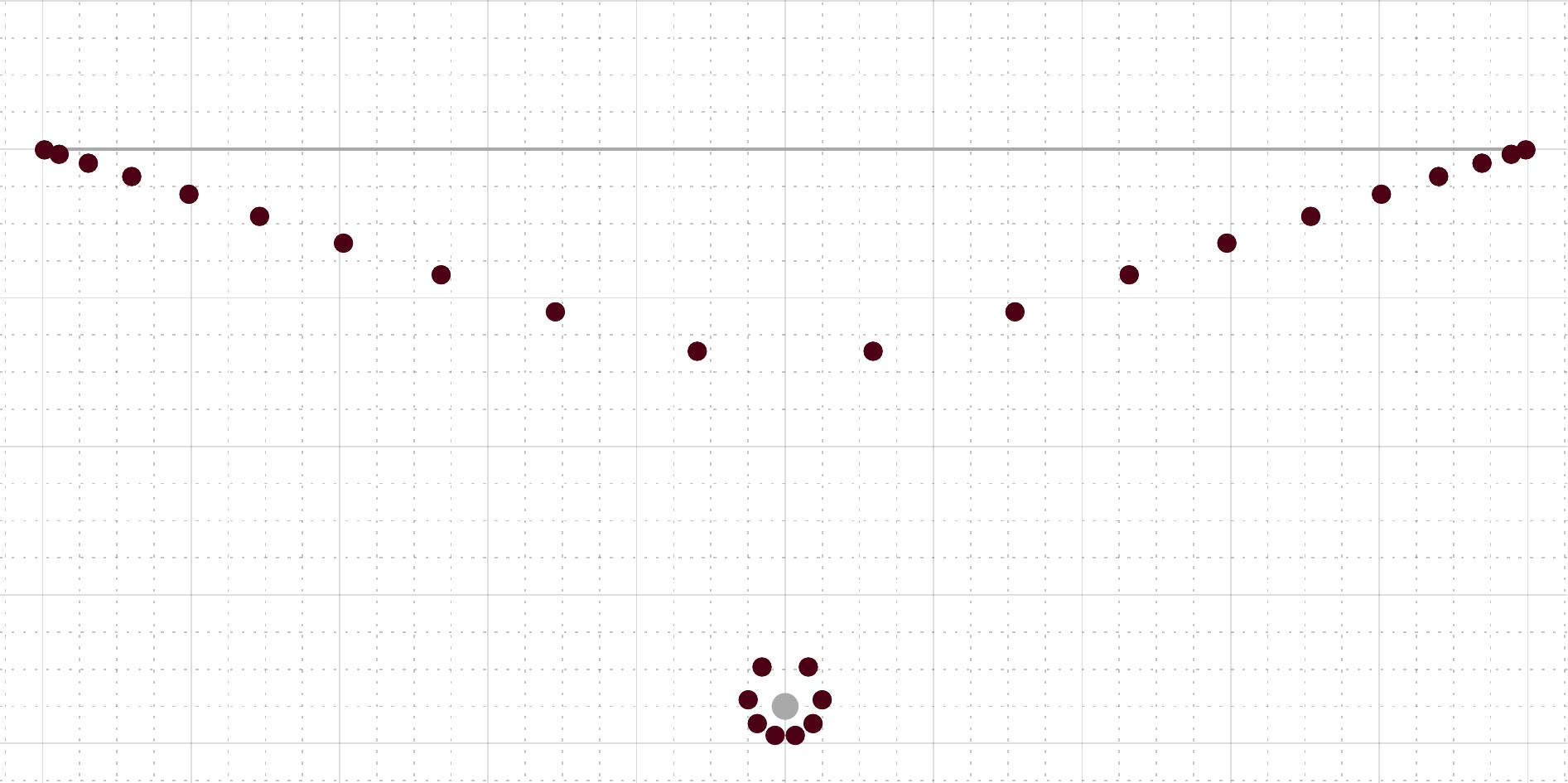}}
\caption{\small Poles of multipoint Pad\'e approximants to \( 1/w_L(z) \) where panel (a): \( L \) is an arc connecting \( -1 \) to some \( x_*>5/4 \) through the upper half-plane and then \( x_* \) to \( 1 \) through the lower half-plane and there are \( 48 \) interpolation conditions at infinity and \( 10 \) conditions at \( 5/4 \); panel (b): \( L \) is a lower unit semi-circle and there are \( 48 \) interpolation conditions at infinity and \( 8 \) conditions at \( -3\ic/4 \).}
\label{fig:1}
\end{figure}

The goal of this note is to explain the numerical results shown on Figure~\ref{fig:1}. To accomplish it, let us set
\[
\varphi_L(z) := z - w_L(z), \quad z\in D_L,
\]
which is an analytic function in \( D_L \) with a simple zero at infinity that satisfies \( \varphi_{L+}(s)\varphi_{L-}(s) \equiv 1 \), \( s\in L \).  Given an interpolation scheme \( \{E_i\}_{i=1}^\infty \), define
\begin{equation}
\label{Bi}
B_i(\zeta) := \prod_{e\in E_i} \frac{\zeta-\varphi_L(e)}{1-\zeta\varphi_L(e)}.
\end{equation}
Clearly, each \( B_i(\zeta) \) is a rational function with \( i \) zeros and \( i \) poles in \( \overline\C \), counting multiplicities. All the zeros of \( B_i(\zeta) \) belong to the interior of a Jordan curve \( J^{-1}(L) \), where  \( J(\zeta) := (\zeta+1/\zeta)/2 \) is the Joukovski map (\( \varphi_L(z) \) maps \( D_L \) conformally onto this domain and \( J(\varphi_L(z))=z \) for \( z\in D_L \)). Notice that if a multi-set \( E_i \) is conjugate-symmetric (i.e., \( e\in E_i \) if and only if \( \overline e\in E_i \)), then \( B_i(\zeta) \) is a Blaschke product. Moreover, if  \( L=[-1,1] \), then all the zeros of \( B_i(\zeta) \), that is, points \( \varphi_{[-1,1]}(e) \), \( e\in E_i \), belong to the unit disk \( \D \).

We shall consider only those interpolation schemes for which the level-lines \( \{|B_i(\tau)|=1\} \) have a definite asymptotic behavior. That is, we assume the following.

\begin{assumption}[colbacktitle=webgreen]{}{1}
Arc \( L \) and interpolation scheme \( \{E_i\}_{i=1}^\infty \) in \( D_L \) are such that
\begin{itemize}
\item[(i)] there exist \( M>0 \) and an integer \( l_*\geq0 \) such that the functions \( B_i(\zeta) \) from \eqref{Bi} satisfy
\[
|B_i(\tau)| \leq M, \quad \tau\in \Gamma := \cup_{l=-l_*}^{l_*}\Gamma_l, \quad i\in\N,
\]
where \( \Gamma_l \)  are pairwise disjoint smooth Jordan curves, \( \Gamma_0 \) contains \( \pm1 \), \( \Gamma_l \), \( l>0 \), belong to the interior of \( \Gamma_0 \), and \( \Gamma_{-l} = \{\tau:\tau^{-1}\in\Gamma_l \} \), \( 0\leq l\leq l_* \);
\item[(ii)] in each connected component of the complement of \( \Gamma \) it holds locally uniformly that either \( |B_i(\zeta)| \to 0 \) as \( i\to\infty \) or \( |B_i(\zeta)| \to \infty \) as \( i\to\infty \).
\end{itemize}
\end{assumption}

It follows from \eqref{Bi} that \( B_i(1/\zeta)=1/B_i(\zeta) \) and therefore \( |B_i(\tau)|=1 \) if and only if \( |B_i(1/\tau)|=1 \). This explains why \( \pm1 \in \Gamma \) and why each connected component of \( \Gamma \) must be either invariant under the map \( \tau\mapsto1/\tau \) (there is only one such component by the assumption, namely, \( \Gamma_0 \)), or is mapped into another connected component by this map. This symmetry also shows that \( |B_i(\tau)| \geq M^{-1}>0 \), \( \tau\in \Gamma \), \( i\in\N \), and of course \( M\geq1 \). 

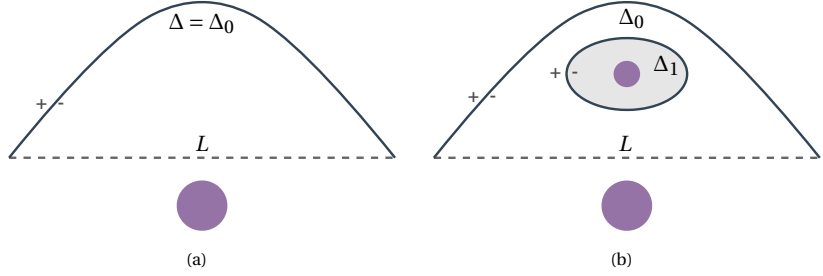
\begin{figure}[ht!]
\centering
\subfigure[]{\includegraphics[scale=1.25]{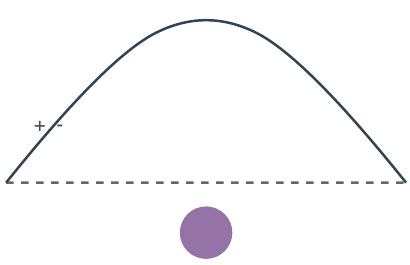}}
\begin{picture}(0,0)
\put(-78,32){\( L \)}
\put(-88,78){\( \Delta=\Delta_0 \)}
\end{picture}
\quad
\subfigure[]{\includegraphics[scale=1.25]{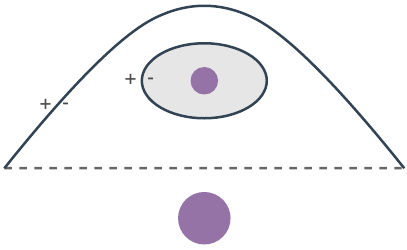}}
\begin{picture}(0,0)
\put(-78,32){\( L \)}
\put(-78,80){\( \Delta_0 \)}
\put(-65,62){\( \Delta_1 \)}
\end{picture}
\caption{\small Darker filled circles represent interpolation points (bigger circle represents more interpolation conditions at the point), dashed lines represent \( L \), solid lines represent \( \Delta \), and lightly shaded regions represent \( D_\Delta^\infty \). Panel (a): interpolation points create an external field that pushes \( L \) up to \( \Delta \). Panel (b): interpolation points below \( L \) push it up, interpolation points above \( L \) push it down, but create weaker external field resulting in \( L \) going through them while simultaneously forming a barrier \( \Delta_1 \).}
\label{fig:2}
\end{figure}

Let \( \Delta := J(\Gamma) \).  We shall say that \( \Delta \) is a \emph{symmetric contour corresponding to the interpolation scheme \( \{E_i\}_{i=1}^\infty \)}. We write \( \Delta = \cup_{l=0}^{l_*} \Delta_l \), where \( \Delta_l := J(\Gamma_l) \). Notice that \( \Delta_0 \) is a Jordan arc connecting \( -1 \) and \( 1 \), and each \( \Delta_l \), \( l>0 \), is a Jordan curve that contains \( \Delta_0 \) in its exterior.  

\begin{figure}[ht!]
\centering
\subfigure[]{\includegraphics[scale=1.25]{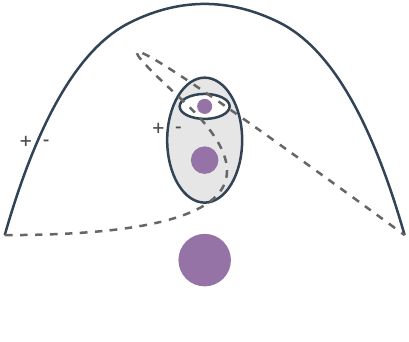}}
\begin{picture}(0,0)
\put(-105,44){\( L \)}
\put(-81,114){\( \Delta_0 \)}
\put(-75,74){\( \Delta_2 \)}
\put(-81,98){\( \Delta_1 \)}
\end{picture}
\quad
\subfigure[]{\includegraphics[scale=1.25]{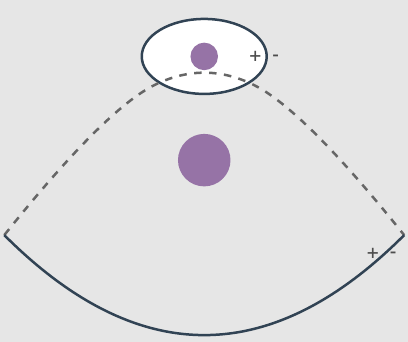}}
\begin{picture}(0,0)
\put(-138,60){\( L \)}
\put(-83,5){\( \Delta_0 \)}
\put(-83,82){\( \Delta_1 \)}
\end{picture}
\caption{\small This is a continuation of Figure~\ref{fig:2}. Panel (a): top and bottom groups of interpolation points create an external field that pushes \( L \) up while the middle group pushes \( L \) down, due to different strength of the components of the external field generated by these groups, two barriers are created. Panel (b): interpolation points below \( L \) create an external field that pushes \( L \) up all the way through \( \infty \) to the displayed position of \( \Delta_0 \), interpolation points above \( L \) create a weaker external field that results in a barrier \( \Delta_1 \).}
\label{fig:3}
\end{figure}
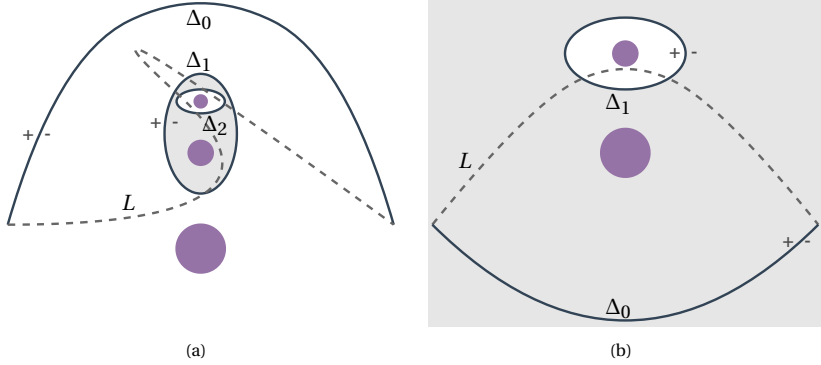

\begin{figure}[ht!]
\centering
\subfigure[]{\includegraphics[scale=0.3]{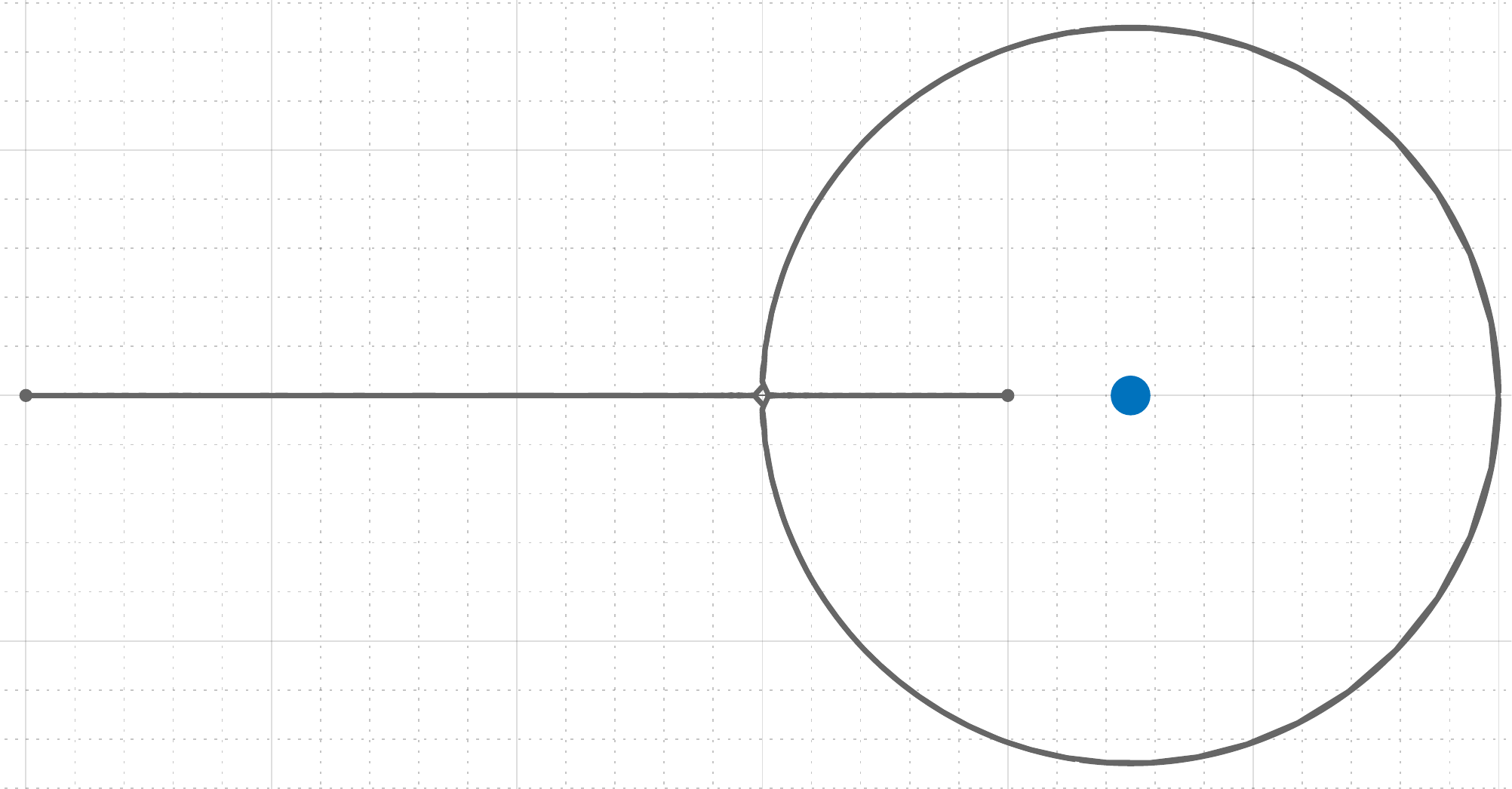}}
\quad
\subfigure[]{\includegraphics[scale=0.3]{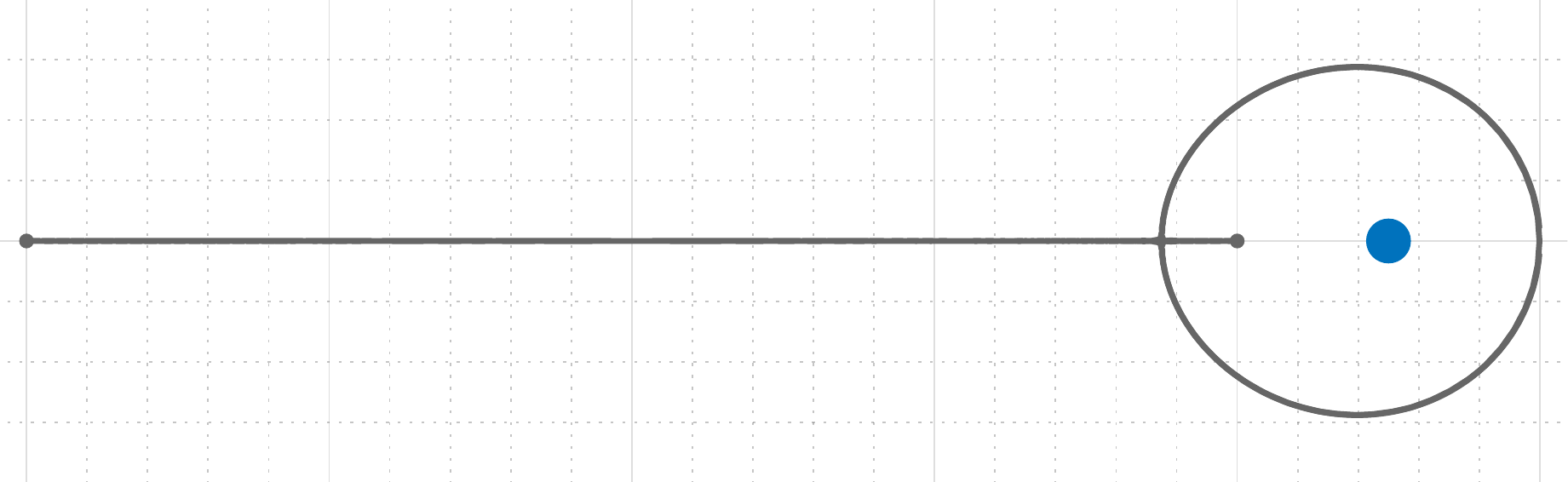}}
\\
\subfigure[]{\includegraphics[scale=0.35]{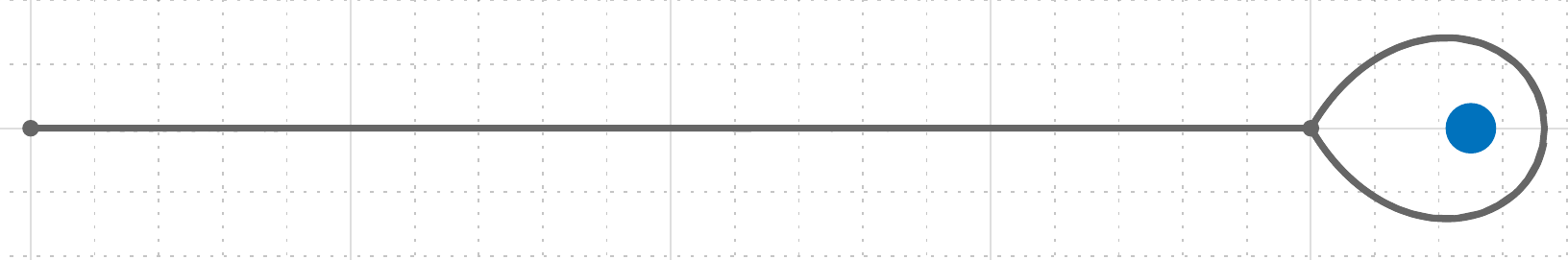}}
\quad
\subfigure[]{\includegraphics[scale=0.35]{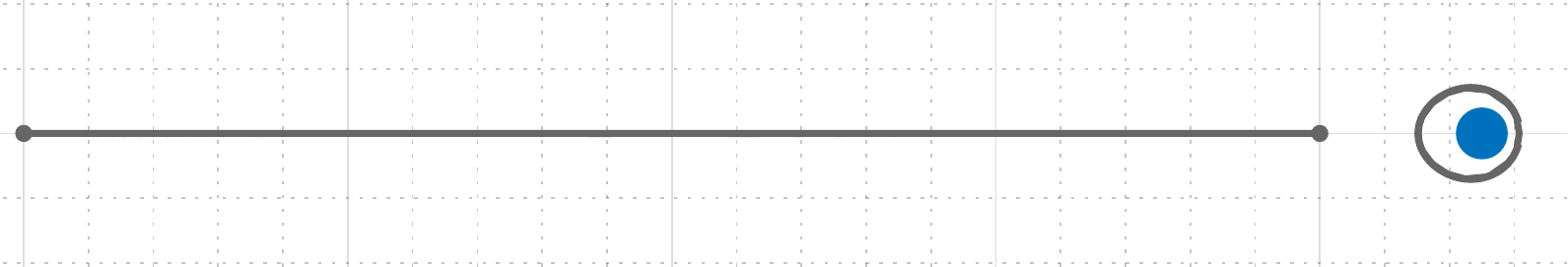}}
\caption{\small Symmetric contours \( \Delta \) that correspond to \( L \) that connects \( -1 \) to some \( x_*>5/4 \) through the upper half-plane and then \( x_* \) to \( 1 \) through the lower half-plane and interpolation schemes where the interpolation conditions are equally distributed  between  between \( \infty \) and \( 5/4 \) (panel (a)) or there are twice (panel (b)), three times (panel (c)), or four times (panel (d)) more interpolation conditions at \( \infty \) than at \( 5/4 \) (the disk on all of the figures). The plots are obtained by numerically plotting the level-line \( \{|B_i(\tau)|=1\} \).}
\label{fig:4}
\end{figure}

\begin{figure}[ht!]
\centering
\subfigure[]{\includegraphics[scale=0.3]{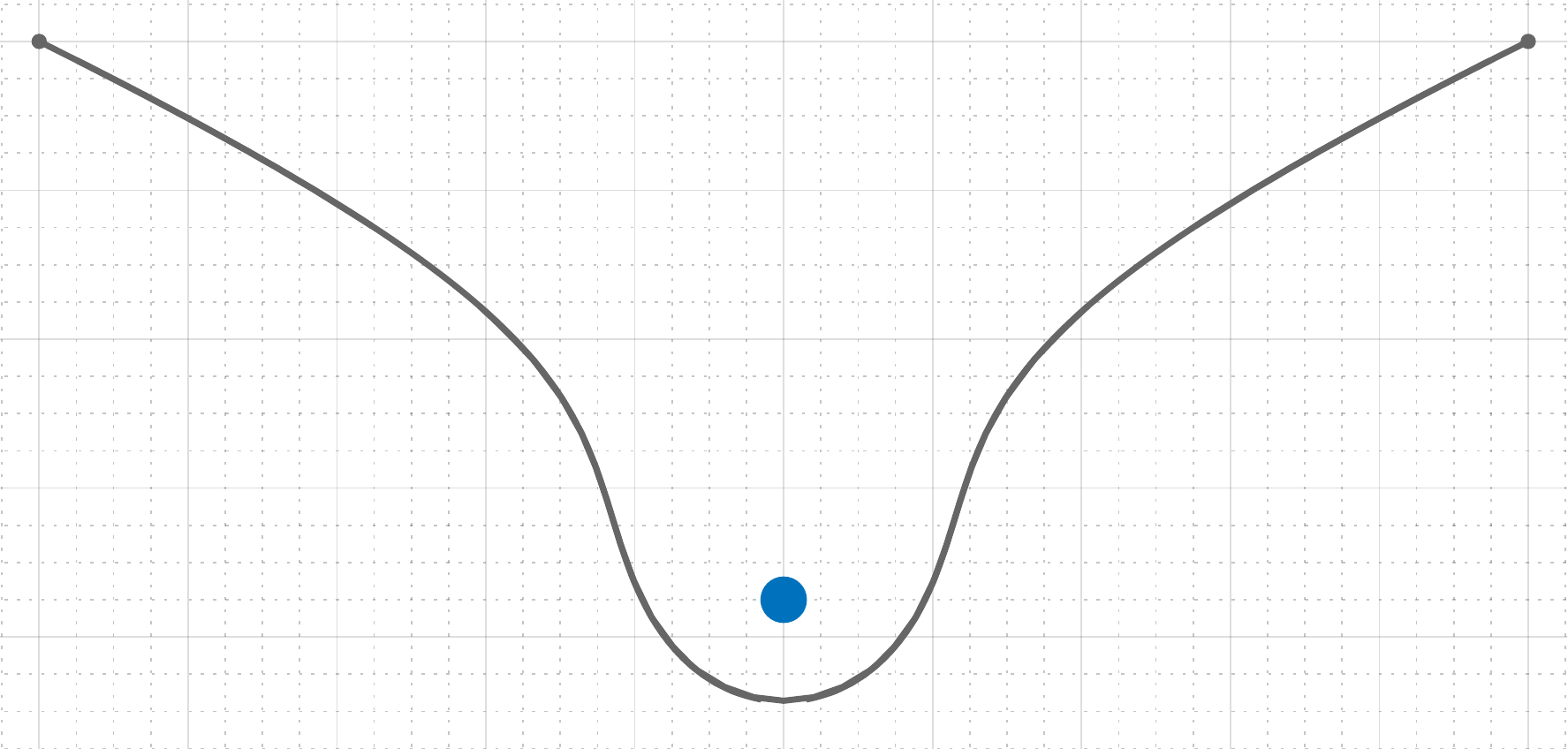}}
\quad
\subfigure[]{\includegraphics[scale=0.3]{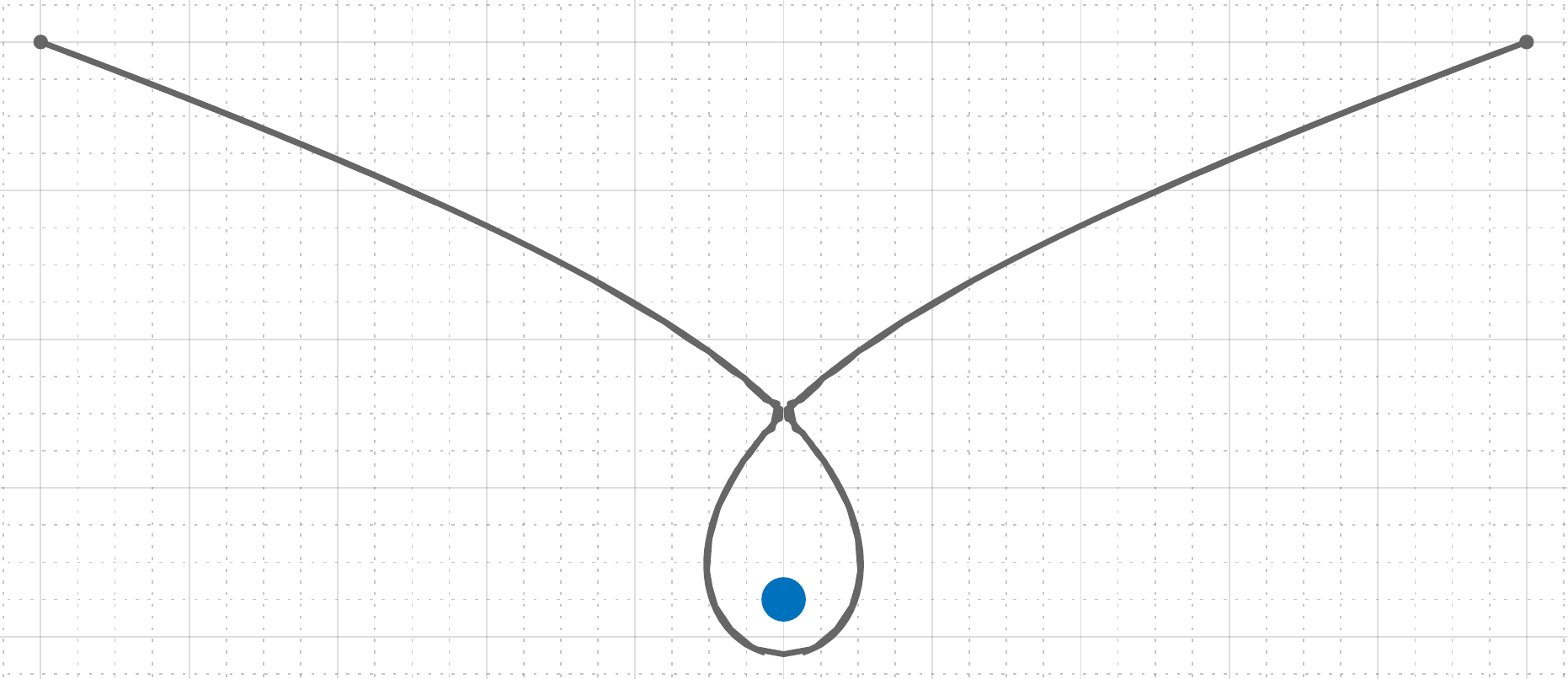}}
\quad
\subfigure[]{\includegraphics[scale=0.35]{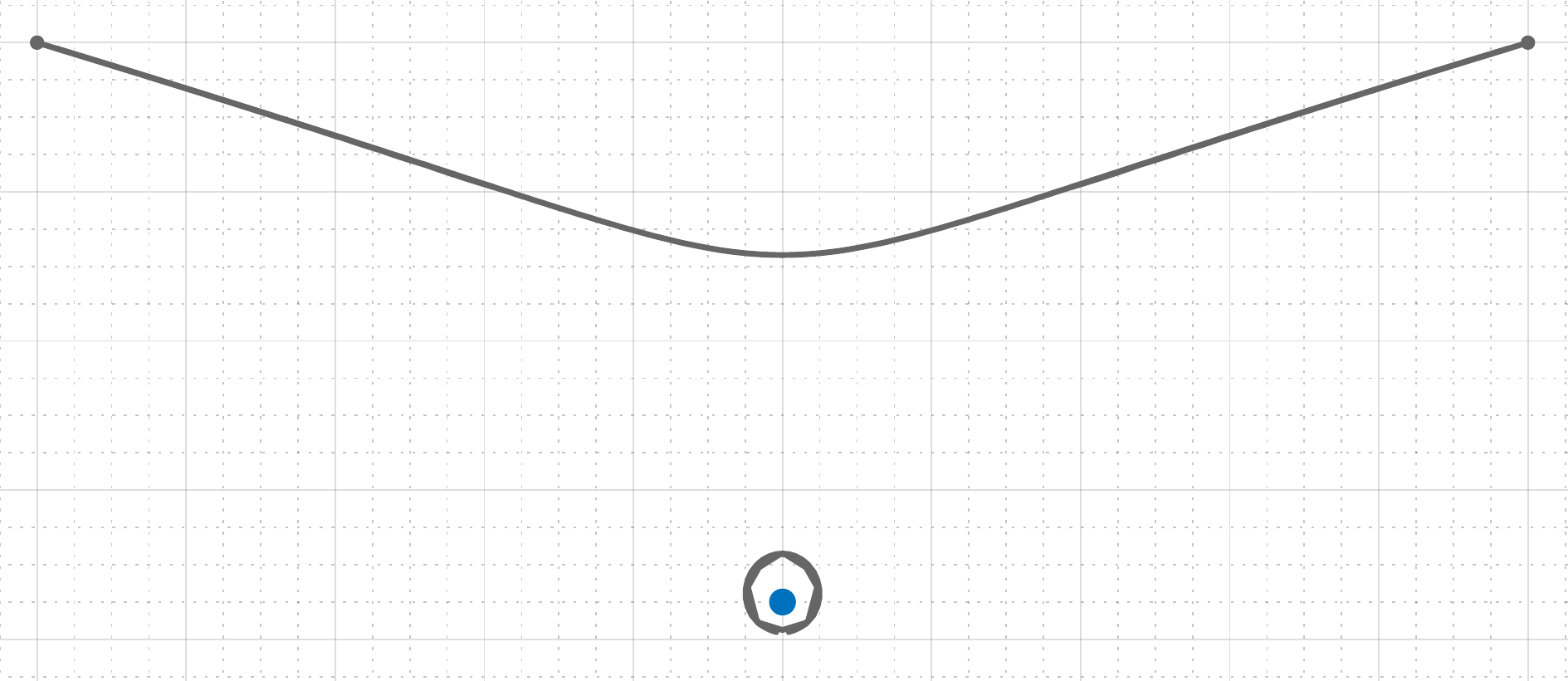}}
\caption{\small Symmetric contours \( \Delta \) that correspond to \( L \) being a lower unit semi-circle and interpolation schemes where there are four (panel (a)), five (panel (b)), or six (panel (c)) times more interpolation conditions at \( \infty \) than at \( -3\ic /4 \) (the disk on all of the figures).}
\label{fig:5}
\end{figure}

One way to think about the symmetric contours defined above is through the level lines of functions \( B_i(\zeta) \) from \eqref{Bi}. Given a contour \( L \), fix \( E_N := \{ e_1,\ldots,e_N\} \), a collection of not necessarily distinct points in \( D_L \). Define \( E_0:=\varnothing \), \( E_i:=\{e_1,\ldots,e_i\} \) for \( i<N \), and \( E_{kN+i}=(kE_N)\cup E_i \), \( k\geq 1 \), \( 0\leq i<N \), where \( kE_N \) is the multi-set consisting of \( k \) copies of \( E_N \). Clearly, it holds that
\[
B_{kN+i}(\zeta) = B_N^k(\zeta) B_i(\zeta)
\]
for any \( k\geq0 \) and \( 0\leq i<N \). Let \( \Gamma := \{ \tau:|B_N(\tau)|=1 \} \) be the \( 1 \)-level line of \( B_N(\zeta) \). As mentioned before, it follows immediately from \eqref{Bi} that \( \Gamma \) is invariant under the map \( \tau \mapsto 1/\tau \). It also follows from the choice of the sets \( E_i \), \( i\leq N \), that there exists \( M\geq 1 \) such that \( M^{-1}\leq |B_i(\tau)|\leq M \) for \( \tau\in\Gamma \) and \( i<N \). Thus, it necessarily holds that \( M^{-1}\leq |B_i(\tau)|\leq M \) for all \( i\in\N \). It is also clear that condition (ii) of Assumption~\ref{as:1} is satisfied. Hence, if \( \Gamma \) consists of \emph{non-intersecting} Jordan curves, it satisfies all the requirements of Assumption~\ref{as:1}. The projection of \( \Gamma \) by the Joukovsky map \( J(\zeta) \) is then an admissible symmetric contour \( \Delta \) that we work with in this note. In general, there are no reasons to exclude the cases of contours \( \Gamma \) with points of self-intersection. We do it here only to simplify the analysis of the asymptotic behavior of the approximants. In particular, Figures~\ref{fig:4} and~\ref{fig:5} depict possible contours \( \Delta \) obtained from the above construction for sets \( E_N \) containing only two distinct element in different proportions (one of them being infinity). In general,  sets \( \Gamma \) from Assumption~\ref{as:1} are not level lines of any particular rational function of the form \eqref{Bi}, but certain limits of such level lines.

Another way to think about the symmetric contours defined above is through heuristics of potential theory. One can think of \( L \) as an elastic band that can freely move in the complex plane except for the fixed points \( \pm1 \). Assume also that, say positive, unit charge can freely float along \( L \). Somewhat unrealistically suppose also that the band can attach to itself to form a loop and then snap into disjoint pieces at the point where the loop was formed. Distribute  a unit, also positive, charge among the points of the set \( E_i \), where each point gets the amount of charge proportional to its multiplicity in \( E_i \). These charges create an external field that makes \( L \) change its location. If \( L \) settles down to some collection of bands in the limit as \( i\to\infty \), this is exactly a symmetric contour \( \Delta \) defined above. This heuristics is the one used in captions to Figures~\ref{fig:2} and~\ref{fig:3}. In particular, if \( L \) is the lower unit semicircle, then the unit charge at infinity will move it to \( \Delta=[-1,1] \). If we place one fifth of the charge at point \( -3\ic/4 \) and keep the rest of it at infinity, then we will get \( \Delta \) as on Figure~\ref{fig:5}(a). By decreasing the charge at \( -3\ic/4 \), we will make \( L \) wrap around \( -3\ic/4 \) until it forms a loop around \( -3\ic/4 \) when \( 1/6 \) the charge is placed there, Figure~\ref{fig:5}(b). If the charge at \( -3\ic/4 \) is further decreased, the band \( L \) snaps into two pieces, one band \( \Delta_0 \) connecting \( -1 \) and \( 1 \) and one loop \( \Delta_1 \) around \( -3\ic/4 \), see Figure~\ref{fig:5}(c).

Let us now continue with the assumption on the weight function \( \rho(s) \) in \eqref{f}. To this end, we introduce an orientation for all the connected components of \( \Delta \). We orient \( \Delta_0 \) from \( -1 \) to \( 1 \). To orient \( \Delta_l \), \( l>0 \), put
\begin{equation}
\label{bi}
b_i(z) := \prod_{e\in E_i} \frac{\varphi_{\Delta_0}(z)-\varphi_L(e)}{1-\varphi_{\Delta_0}(z)\varphi_L(e)} = B_i(\zeta),
\end{equation}
where \( z=J(\zeta) \) and \( \zeta \) belongs to the interior domain of \( \Gamma_0 \). Let \( D_\Delta := \overline\C\setminus\Delta \) and denote by \( D_\Delta^0 \) and \( D_\Delta^\infty \) the parts of \( D_\Delta \) where the functions \( |b_i(z)| \) converge to \( 0 \) and diverge to \( \infty \), respectively (\( D_\Delta = D_\Delta^0 \cup D_\Delta^\infty \)), see Figures~\ref{fig:2} and~\ref{fig:3}.  Notice that each curve \( \Delta_l \), \( l>0 \), is a part of the boundary of both \( D_\Delta^0 \) and \( D_\Delta^\infty \) and we orient it so that \( D_\Delta^0 \) remains on the left as \( \Delta_l \) is traversed in the positive direction (this orientation allows us to distinguish \( + \) side (left) and \( - \) side (right) of \( \Delta_l \) in a manner convenient for our purposes; our choice is not related to the location of the interior domain of \( \Delta_l \), it might end up on any side of the curve, see, for example, Figure~\ref{fig:3}). Let us write \( \overline\C\setminus(L\cup\Delta_0) =: U_u \cup U_b \), where \( U_u \) is the unbounded component of \( \overline\C\setminus(L\cup\Delta_0) \).  We shall also set \( U_{l,b} \) to be the interior of \( \Delta_l \), \( l>0 \).

\begin{assumption}[colbacktitle=webgreen]{}{2}
It is assumed that \( \rho(z) \) is analytic in a domain that contains each \( \overline U_{l,b} \), \( l>0 \), as well as \( \overline U_b \) (resp. \( \overline U_u \)) when \(\infty \in D_\Delta^0 \) (resp. \( \infty\in D_\Delta^\infty \)). It is also assumed that \( \rho(z) \) is non-vanishing on \( \Delta_0 \) and \( \cup_{l>0}\overline U_{l,b} \) (thus, it has zero winding number on each \( \Delta_l \), \( l>0 \)).
\end{assumption}

The situation \(\infty \in D_\Delta^0 \) is schematically depicted on Figures~\ref{fig:2}(a,b) and~\ref{fig:3}(a). It means that the unbounded component of the complement of \( D_\Delta \) is a part of \( D_\Delta^0 \). The situation \(\infty \in D_\Delta^\infty \) is schematically depicted on Figure~\ref{fig:3}(b) and represents geometries where the unbounded component of the complement of \( D_\Delta \) is a part of \( D_\Delta^\infty \). Under Assumption~\ref{as:2} the following proposition holds.

\begin{proposition}[colbacktitle=RoyalBlue]{}{1}
For all \( i \) large enough each element of \( E_i \) has a neighborhood in which \( \widehat\rho_L(z)=\widehat\rho_\Delta(z) \), where
\begin{equation}
\label{rhoDelta}
\widehat\rho_\Delta(z) := \frac1{2\pi\ic}\int_\Delta\frac{\rho(s)}{s-z}\frac{\dd s}{w(s)}, \quad z\in D_\Delta,
\end{equation}
we set \( \varsigma := 1 \) when \( \infty\in D_\Delta^0 \) and \( \varsigma := -1 \) when \( \infty\in D_\Delta^\infty \), as well as
\begin{equation}
\label{w}
w(s) := \left\{
\begin{array}{rl}
\varsigma w_{\Delta_0+}(s), & s\in\Delta_0, \medskip \\
w_{\Delta_0}(s), & s\in\cup_{l>0}\Delta_l.
\end{array}
\right.
\end{equation}
In particular, \( [m/n;E_{m+n}]_{\widehat\rho_L}(z)=[m/n;E_{m+n}]_{\widehat\rho_\Delta}(z) \) for all \( m+n \) large enough.
\end{proposition}

In view of the above proposition we can now define functions of the second kind as
\begin{equation}
\label{SKind}
R_{m,n}(z) := \frac{q_{m,n}(z)\widehat\rho_\Delta(z)-p_{m,n}(z)}{v_{m,n}(z)}, \quad z\in D_\Delta
\end{equation}
(in the diagonal case we shall simply denote them by \( R_n(z) \)). It follows from \eqref{Rn} and Proposition~\ref{prop:1} that \( R_{m,n}(z) \) is analytic in the domain of its definition and vanishes at infinity with order at least \( \min\{m,n\}+1 \).

\begin{proposition}[colbacktitle=RoyalBlue]{}{2}
With \( w(s) \) given by \eqref{w}, it holds that
\[
R_{m,n}(z) = \frac1{2\pi\ic}\int_\Delta\frac{q_{m,n}(s)}{v_{m,n}(s)}\frac{\rho(s)}{w(s)}\frac{\dd s}{s-z}, \quad z\in D_\Delta.
\]
\end{proposition}

Recall that \( \rho(s) \) is assumed to be non-vanishing on \( \Delta \) with a zero winding number on each \( \Delta_l \), \( l>0 \). Thus, we can define \( \log\rho(s) \) continuously on \( \Delta \). Let 
\begin{equation}
\label{Sz}
S_\rho(z) := \exp\left\{\frac{w_{\Delta_0}(z)}{2\pi\ic}\int_\Delta\frac{\log\rho(s)}{z-s}\frac{\dd s}{w(s)}\right\}, \quad z\in D_\Delta,
\end{equation}
be the Szeg\H{o} function of \( \rho(s) \). That is, \( S_\rho(z) \) is holomorphic and non-vanishing in \( D_\Delta \) with continuous traces on \( \Delta \) that satisfy
\begin{equation}
\label{Sz-j}
\left\{
\begin{array}{ll}
S_{\rho+}(s)S_{\rho-}(s) = \rho^{-\varsigma}(s), & s\in \Delta_0, \medskip \\
S_{\rho+}(s) = S_{\rho-}(s)/\rho(s), & s\in \Delta_l, ~~l>0.
\end{array}
\right.
\end{equation}

\begin{proposition}[colbacktitle=RoyalBlue]{}{3}
Let \( b_{2n}(z) \) be given by \eqref{bi} and \( S_\rho(z) \) be given by \eqref{Sz}. Then there exist choices of the branch of the square root so that
\begin{equation}
\label{Psin}
\Psi_n(z) := 
\left\{
\begin{array}{rl}
\sqrt{v_n(z)/b_{2n}(z)}S_\rho(z), & z\in D_\Delta^0, \medskip \\
\sqrt{v_n(z)b_{2n}(z)}/S_\rho(z), & z\in D_\Delta^\infty,
\end{array}
\right.
\end{equation}
is an analytic and non-vanishing function in \( D_\Delta \) that satisfies
\begin{equation}
\label{P-pr}
\left\{
\begin{array}{ll}
\Psi_n(z) = \gamma_n^{-1}z^n + \mathcal O(z^{n-1}) & \text{as} \quad z\to\infty, \medskip \\
\Psi_{n+}(s)\Psi_{n-}(s) = (v_n/\rho)(s), & s\in \Delta,
\end{array}
\right.
\end{equation}
for some non-zero constant \( \gamma_n \). 
\end{proposition}

The following theorem constitutes the main result of this note.

\begin{theorem}[colbacktitle=webbrown]{}{1}
Let \( L \), \( \{E_i\}_{i=1}^\infty \), and \( \rho(z) \) be such that Assumptions~\ref{as:1} and~\ref{as:2} are satisfied. Further, let \( \Delta \) be the symmetric contour corresponding to the interpolation scheme \( \{E_i\}_{i=1}^\infty \) and \( p_n(z)/q_n(z)=[n/n;E_{2n}]_{\widehat\rho_L}(z) \) be the sequence of diagonal multipoint Pad\'e approximants of \( \widehat\rho_L(z) \) associated with \(  \{E_i\}_{i=1}^\infty \), see \eqref{Rn}. Then it holds uniformly on closed subsets of \( D_\Delta \) that
\begin{equation}
\label{asymp1}
q_n(z) = \gamma_n (1+o(1))\Psi_n(z)
\end{equation}
and
\begin{equation}
\label{asymp2}
R_n(z) = \frac{\gamma_n(1+o(1))}{w_{\Delta_0}(z)\Psi_n(z)}
\left\{
\begin{array}{rl}
1, & z\in D_\Delta^0, \medskip \\
-1, & z\in D_\Delta^\infty.
\end{array}
\right.
\end{equation}
In particular, it holds uniformly on closed subsets of \( D_\Delta \) that
\begin{equation}
\label{asymp3}
\widehat\rho_\Delta(z) - [n/n;E_{2n}]_{\widehat\rho_L}(z) = \frac{1+o(1)}{w_{\Delta_0}(z)}
\left\{
\begin{array}{rl}
b_{2n}(z)S_\rho^{-2}(z), & z\in D_\Delta^0, \medskip \\
-b_{2n}^{-1}(z)S_\rho^2(z), & z\in D_\Delta^\infty.
\end{array}
\right.
\end{equation}
\end{theorem}

Let us place the above theorem, which covers contours on Figures~\ref{fig:4}(d) and \ref{fig:5}(c), in the context of recent and not so distant results on attracting curves of the poles of Pad\'e approximants. When the approximants  are classical (i.e., all interpolation points are at infinity) and the approximated functions are of hyperelliptic type (Cauchy integrals of densities of the form \( \rho(s)/w_{L+}(s) \), where \( w_L^2(s) \) is a polynomial of even degree), it was recognized by Nuttall and Singh \cite{NutS77} that the pole attracting contour \( \Delta \) is a projection of a certain level line on the Riemann surface of \( w_L(z) \) to the complex plane (such \( \Delta \)'s do not separate the plane). In this case an analog of Theorem~\ref{thm:1} was obtained in \cite{Ya15} by the author (see also \cite{Suet00,Suet02,BYa13}). This approach was later applied to multipoint Pad\'e approximants by Baratchart and the author in \cite{BYa09c} when \(\Delta\) is Jordan arc, which covers the contour on Figure~\ref{fig:5}(a) (see also \cite{CYLL98,St99})  and by the author \cite{Ya18} in the full generality of \cite{NutS77} for contours that do not separate the plane. So, this work is a natural continuation of  \cite{BYa09c} (\(\Delta_0\) is a Jordan arc in Theorem~\ref{thm:1}) and the author is currently preparing generalization of \cite{Ya18} to the case of contours that do disconnect the plane. When the approximated function is multi-valued but not hyperelliptic, the approach of \cite{NutS77} no longer works. However, in a series of pathbreaking papers \cite{St85, St85b, St86} Stahl showed that attracting contours for classical Pad\'e approximants can be identified as the ones having minimal logarithmic capacity among all possible branch cuts for the approximated function. Strong asymptotics of classical Pad\'e approximants for "generic" geometries then was obtained by Aptekarev and the author \cite{ApYa15}. Stahl's approach was extended to the set up that includes multipoint Pad\'e approximants by Gonchar and Rakhmanov \cite{GRakh87} (contours are not allowed to separate the plane). Further extension of Stahl's ideas to contours that are allowed to separate the plane was obtained by Buslaev \cite{Bus13,Bus15}. Buslaev's existence theorem covers all possible contours arising from interpolation schemes with finitely many points (in particular, it covers Figures~\ref{fig:4} and~\ref{fig:5}). However, his convergence results consider only those geometries where each component of \( \C\setminus\Delta \) contains a branch point of the approximated function as on Figure~\ref{fig:4}(a,b) (they are also concerned with weak asymptotics). Some strong asymptotics analogs of Buslaev convergence results were recently obtained by the author \cite{uYa}. In this sense the current note is a continuation of the work done in \cite{uYa}.

\section{Proofs}

Let us recall one of the Plemelj-Sokhotski formulae \cite[Section~I.4.2]{Gakhov} that states that
\[
\Phi_+(s) - \Phi_-(s) = \phi(s), \quad s\in I^\circ, \quad \Phi(z):=\frac1{2\pi\ic}\int_I \frac{\phi(s)\dd s}{s-z},
\]
where \( I \) is a smooth oriented Jordan arc (closed or open), \( I^\circ \) is its interior, and \( \phi(s) \) is a function that is H\"older smooth  on any subarc of \( I^\circ \). It also holds that \( \Phi_\pm(s) \) are H\"older smooth functions on any subarc of \( I^\circ \). Recall further Yamashita-Bagemihl analytic continuation principle, see  \cite[Theorem~3]{Yama69}: if domains \( D_1,D_2\) are disjoint, \( I \) is an open analytic arc that belongs to the boundary of both domains, \( f_i(z) \) belongs to the Hardy space \( H_1(D_i)\), \( i\in\{1,2\} \) (in particular, bounded in \( D_i \)), and
\[
\lim_{ L_1(s)\ni z\to s} f_1(z) = \lim_{ L_2(s)\ni z\to s} f_2(z) \neq \infty
\]
for almost every \( s\in  I \), where \( L_i(s) \) is some arc in \( D_i \) terminating at \( s \), then \( f_2(z) \) is the analytic continuation of \( f_1(z) \) across \( L \) from \( D_1 \) into \( D_2 \).

\begin{proof}[Proof of Proposition~\ref{prop:1}]

Denote the unbounded component of \( D_\Delta \) by \( \Omega \). We shall consider only those indices \( i \) that are large enough so that \( b_i(e)=0 \) when \( e\in E_i\cap D_\Delta^0 \) and \( b_i(e)=\infty\) when \( e\in E_i\cap D_\Delta^\infty \). Fix such an index \( i \). Since \( E_i \) is finite and \( L \) can be perturbed in \eqref{f} without changing the relevant values of \( \widehat\rho_L(z) \) in \eqref{Rn}, we can assume that \( L \) and \( \Delta_0 \) have at most finitely many points in common and these arcs intersect at those points transversally.

According to Plemelj-Sokhotski formulae it follows from \eqref{f} that
\[
\widehat\rho_{L+}(s) - \widehat\rho_{L-}(s) = \rho(s)/w_{L+}(s), \quad s\in L.
\]
Therefore, the analytic continuation of \( \widehat\rho_L(z) \) across \( L \) is given by \( \widehat\rho_L(z) - \rho(z)/w_L(z) \). Define \( \widehat\rho_{\Delta_0}(z) \) as in \eqref{f} with \( L \) replaced by \( \Delta_0 \). Note that \( w_L(z)=-w_{\Delta_0}(z) \) for \( z\in U_b \) and \( w_L(z)=w_{\Delta_0}(z) \) for \( z\in U_u \). Respectively, it also holds that \( \varphi_{\Delta_0}(z) = 1/\varphi_L(z) \) for \( z\in U_b \) and \( \varphi_{\Delta_0}(z) = \varphi_L(z) \) for \( z\in U_u \).

Assume that \( \Omega\subseteq D_\Delta^0 \).   Then
\begin{equation}
\label{1-0}
\widehat\rho_{\Delta_0}(z) = \left\{
\begin{array}{ll}
\widehat\rho_L(z) + \rho(z)/w_{\Delta_0}(z), & z\in U_b, \smallskip \\
\widehat\rho_L(z), & z\in U_u.
\end{array}
\right.
\end{equation}
Indeed, the right- and left-hand sides of \eqref{1-0} are analytic in \( \overline\C\setminus\Delta_0 \), vanish at infinity, and have smooth traces on \( \Delta_0 \).  The jump of the right-hand side across \( \Delta_0 \) is equal to \( \rho(s)/w_{\Delta_0+}(s) \) when \( U_b \) lies to the left of \( \Delta_0 \) and again \( -\rho(s)/w_{\Delta_0-}(s)=\rho(s)/w_{\Delta_0+}(s) \) when \( U_b \) lies to the right of \( \Delta_0 \). That is, in any case it matches the jump of \( \widehat\rho_{\Delta_0}(z) \) across \( \Delta_0 \). Hence, the difference of the left- and right-hand sides of \eqref{1-0} is analytic in \( \C\setminus\{\pm1 \} \) by the analytic continuation principle. As it vanishes at infinity and can have at most square root singularities at \( \pm1 \) by \cite[Section~I.8.4]{Gakhov}, it must be identically zero.  Let \( e\in E_i\cap \Omega \). As \( \Omega\subseteq D_\Delta^0 \), \( b_i(e)=0 \) and therefore \( e\in U_u \). In this case
\begin{equation}
\label{1-1}
\widehat\rho_\Delta(z) = \widehat\rho_{\Delta_0}(z) = \widehat\rho_L(z)
\end{equation}
for \( z \) around \( e \) as claimed, where the first conclusion follows from Assumption~\ref{as:2} and Cauchy theorem. Let now \( e\in  E_i\cap D_\Delta^0 \), \( e\not\in \Omega \). Again, it must hold that \( e\in U_u \). Since each curve \( \Delta_l \), \( l>0 \), lies on the border of both \( D_\Delta^0 \) and \( D_\Delta^\infty \), \( e \) is an interior point of even number of such curves. As orientations of these curves alternate between clockwise and counter-clockwise and 
\begin{equation}
\label{1-2}
\widehat\rho_{\Delta_l}(z) := \frac1{2\pi\ic}\int_{\Delta_l}\frac{\rho(s)}{s-z}\frac{\dd s}{w_{\Delta_0}(s)} = \pm \frac{\rho(z)}{w_{\Delta_0}(z)}, \quad z\in U_{l,b}, 
\end{equation}
(the choice of the sign depends on the orientation of \( \Delta_l \)) by Assumption~\ref{as:2} and Cauchy integral formula, we again get that \eqref{1-1} holds. Now, when \( e\in  E_i\cap D_\Delta^\infty \), \( b_i(z) \) must have a pole at \( e \), that is, \( \varphi_{\Delta_0}(e) = 1/\varphi_L(e) \), which implies that \( e\in U_b \). Let \( \Delta_{l_e} \) be the curve that borders \( \Omega \) and contains \( e \) in its interior. Then \( \Delta_{l_e} \) is clockwise oriented. There could be other curves \( \Delta_l \) containing \( e \) in their interior, but they are even in number and their net contribution to \( \widehat\rho_\Delta(z) \) around \( e \) is zero due to alternating orientations and \eqref{1-2}. Hence,
\[
\widehat\rho_\Delta(z) = \widehat\rho_{\Delta_0}(z) + \widehat\rho_{\Delta_{l_e}}(z) = \widehat\rho_L(z) + \rho(z)/w_{\Delta_0}(z) + \widehat\rho_{\Delta_{l_e}}(z) = \widehat\rho_L(z)
\]
for \( z \) around \( e \) as claimed, where we used \eqref{1-0} and \eqref{1-2} one more time.

Assume now that \( \Omega\subseteq D_\Delta^\infty \). Since \( w_L(z) =w_{\Delta_0}(z) \) for \( z\in U_u \), reasoning as in \eqref{1-0} now implies that
\begin{equation}
\label{1-3}
-\widehat\rho_{\Delta_0}(z) = \left\{
\begin{array}{ll}
\widehat\rho_L(z), & z\in U_b, \smallskip \\
\widehat\rho_L(z) - \rho(z)/w_{\Delta_0}(z), & z\in U_u.
\end{array}
\right.
\end{equation}
Let \( e\in E_i\cap \Omega \). Since \( \Omega\subseteq D_\Delta^\infty \), \( b_i(z) \) has a pole at \( e \) and therefore \( e\in U_b \). Recall \eqref{w}. Since \( \rho(z) \) is analytic on each \( \overline U_{l,b} \), we get as in \eqref{1-1} that
\begin{equation}
\label{1-4}
\widehat\rho_\Delta(z) = -\widehat\rho_{\Delta_0}(z) = \widehat\rho_L(z)
\end{equation}
for \( z \) around \( e \). If \( e\in E_i\cap ( D_\Delta^\infty\setminus\Omega) \), then \( e\in U_b \) and it lies interior to even number of curves \( \Delta_l \), \( l>0 \),  with consecutively alternating orientations. Then the net contribution to \( \widehat\rho_\Delta(z) \) by these curves is zero and \eqref{1-4} remains valid. Finally, if \( e\in E_i\cap D_\Delta^0 \), then \( e\in U_u \) and it lies interior to odd number of curves \( \Delta_l \), \( l>0 \),  with consecutively alternating orientations starting (and ending) with the counter-clockwise one. The desired claim now follows from \eqref{1-2} and \eqref{1-3}.
\end{proof}

\begin{proof}[Proof of Proposition~\ref{prop:2}]
It follows from \eqref{rhoDelta}, \eqref{SKind}, Plemelj-Sokhotski formulae, and the analytic continuation principle that the difference
\[
R_{m,n}(z) - \frac1{2\pi\ic}\int_\Delta\frac{q_{m,n}(s)}{v_{m,n}(s)}\frac{\rho(s)}{w(s)}\frac{\dd s}{s-z}
\]
is analytic in \( \overline\C\setminus\{\pm1\} \) and is vanishing at infinity. It is also follows from the known behavior of Cauchy integrals around the endpoints of contours of integration \cite[Section~I.8.4]{Gakhov} that it can have at most square root singularity at \( \pm1 \). Hence, it must be identically zero.
\end{proof}

\begin{proof}[Proof of Proposition~\ref{prop:3}]
Write \( v_n(z) = v_{n,0}(z)v_{n,\infty}(z) \), where \( v_{n,0}(z) \) and \( v_{n,\infty}(z) \) are monic polynomials that vanish at the zeros and poles of \( b_{2n}(z) \) in \( \C \), respectively. It follows from \eqref{bi} that
\[
b_{2n}(z)\frac{v_{n,\infty}(z)}{v_{n,0}(z)} = \left( c_n \zeta^{n-\deg(v_{n,\infty})}\frac{\prod_{e\in E_{2n,\infty}}(\zeta-\varphi_L(e))}{\prod_{e\in E_{2n,0}}(1-\zeta\varphi_L(e))} \right)^2 := T_n^2(\zeta),
\]
where \( z=J(\zeta) \), \( \zeta \) belongs to the interior domain of \( \Gamma_0 \), \( E_{2n,0} \) and \( E_{2n,\infty} \) are the multi-sets of zeros of \( v_{n,0}(z) \) and \( v_{n,\infty}(z) \) (i.e., counting multiplicities), and
\[
 c_n^2 := (-2)^{\deg(v_{n,0})-\deg(v_{n,\infty})}\prod_{e\in E_{2n,0}}\varphi_L(e)\prod_{e\in E_{2n,\infty}}\varphi_L^{-1}(e)
\]
whose square root we fix arbitrarily. Observe that \( T_n(\zeta) \) is a non-vanishing holomorphic function in the interior of \( \Gamma_0 \) except for a possible zero or pole at the origin. In any case, we can define
\[
\sqrt{(b_{2n}v_{n,\infty}/v_{n,0})(z)} := T_n(\zeta),
\]
where \( z=J(\zeta) \) and \( \zeta \) belongs to the interior domain of \( \Gamma_0 \), which is a holomorphic and non-vanishing function in \( \C\setminus\Delta_0 \) that has either a pole or a zero at infinity. Then we have that
\begin{equation}
\label{3-1}
\sqrt{(b_{2n}v_{n,\infty}/v_{n,0})(s)}_+\sqrt{(b_{2n}v_{n,\infty}/v_{n,0})(s)}_- = T_n(\tau) T_n(1/\tau) =  (v_{n,\infty}/v_{n,0})(s)
\end{equation}
for \( s= J(\tau)\in \Delta_0 \), i.e., \( \tau\in\Gamma_0 \). Now we can define
\begin{equation}
\label{3-2}
\sqrt{v_n(z)/b_{2n}(z)} := v_{n,\infty}(z)T_n^{-1}(\zeta) \qandq \sqrt{v_n(z)b_{2n}(z)} := v_{n,0}(z)T_n(\zeta),
\end{equation}
where \( z=J(\zeta) \) and \( \zeta \) belongs to the interior domain of \( \Gamma_0 \). Relations \eqref{P-pr} now easily follow from \eqref{Sz-j},  \eqref{3-1}, and \eqref{3-2} (observe that it necessarily holds that \( E_{2n}=E_{2n}^* \) when \( \infty\in D_\Delta^{\infty} \), in which case \( 2n=\deg(v_n)=\deg(v_{n,0})+\deg(v_{n,\infty})\)).
\end{proof}

\begin{proof}[Proof of Theorem~\ref{thm:1}]

For brevity, let us set
\[
\boldsymbol I := \left(\begin{matrix} 1 & 0 \\ 0 & 1 \end{matrix}\right) \qandq \sigma_3 := \left(\begin{matrix} 1 & 0 \\ 0 & -1 \end{matrix}\right).
\]
To prove the theorem, we follow by now classical approach of Fokas, Its, and Kitaev \cite{FIK91,FIK92} connecting orthogonal polynomials to matrix Riemann-Hilbert problems and then utilizing the non-linear steepest descent method of Deift and Zhou \cite{DZ93}. 

{\bf Step 1.} Consider the following \( 2\times2 \) Riemann-Hilbert problem (\rhy):
\begin{myitemize}
\label{rhy}
\item[(a)] \( \boldsymbol Y(z) \) is analytic in \( \overline\C\setminus\Delta \) and \( \displaystyle \lim_{z\to\infty} \boldsymbol Y(z)z^{-n\sigma_3} = \boldsymbol I \);
\item[(b)] \( \boldsymbol Y(z) \) has continuous traces on \( \Delta\setminus\{\pm1\} \) that satisfy
\[
\displaystyle \boldsymbol Y_+(s) = \boldsymbol Y_-(s) \left(\begin{matrix} 1 & (\rho/(v_nw))(s) \medskip \\ 0 & 1 \end{matrix}\right), \quad s\in \Delta\setminus\{\pm1\},
\]
where \( w(s) \) was defined in \eqref{w};
\item[(c)] it holds that \( \boldsymbol Y(z) = \mathcal O\left(\begin{matrix} 1 & |z-e|^{-1/2} \smallskip \\ 1 & |z-e|^{-1/2}\end{matrix}\right) \) as \( D_\Delta\ni z\to e\in\{\pm1\} \).
\end{myitemize}
To connect \hyperref[rhy]{\rhy} to the polynomials \( q_n(z) \), recall \eqref{SKind}. If it holds that
\begin{equation}
\label{assumption}
\deg(q_n)=n \qandq R_{n+1,n-1}(z)=k_n^{-1}z^{-n}[1+o(1)]  \qasq z\to\infty
\end{equation}
for a non-zero finite constant \( k_n \), then \hyperref[rhy]{\rhy} is solved by
\begin{equation}
\label{Y}
\boldsymbol Y(z) = \left(\begin{matrix}
q_n(z) & R_n(z) \medskip \\
k_n q_{n+1,n-1}(z) & k_nR_{n+1,n-1}(z)
\end{matrix}\right).
\end{equation}
Conversely, if \hyperref[rhy]{\rhy} is solvable, then its solution necessarily has form \eqref{Y} and \eqref{assumption} is satisfied. Indeed, if \eqref{assumption} holds, then so is \hyperref[rhy]{\rhy}(a) by \eqref{Rn} and since \( \deg(q_{n+1,n-1})<n \). The fact that \hyperref[rhy]{\rhy}(b) is fulfilled follows from Proposition~\ref{prop:2} and Plemelj-Sokhotski formulae. Finally, \hyperref[rhy]{\rhy}(c) holds due to the known behavior of Cauchy integrals around endpoints of contours of integration \cite[Section~I.8.4]{Gakhov}. Conversely, if \( \boldsymbol Y(z) \) is a solution of \hyperref[rhy]{\rhy}, then its unique (if \( \boldsymbol Y_1(z) \) and  \( \boldsymbol Y_2(z) \) are solutions, then their determinants are identically equal to \( 1 \) and \( \boldsymbol Y_1(z)\boldsymbol Y_2^{-1}(z) \) must be entire and equal to \( \boldsymbol I \) at infinity, i.e., equal to \( \boldsymbol I \) everywhere). Furthermore, \( q(z) \), the \((1,1)\)-entry of \( \boldsymbol Y(z) \), must be a monic polynomials of degree exactly \( n \) by \hyperref[rhy]{\rhy}(a,b) and \( R(z) \), the \((1,2)\)-entry, must have an integral representation as in Proposition~\ref{prop:2} with \( q_n(z) \) replaced by \( q(z) \) by \hyperref[rhy]{\rhy}(b,c) and the analytic continuation principle. Replace the first row of \( \boldsymbol Y(z) \) with
\[
\left(\begin{matrix}  q(z)+q_n(z) & R(z) + R_n(z) \end{matrix}\right) \quad \text{or} \quad \left(\begin{matrix}  (q(z)+q_n(z))/2 & (R(z) + R_n(z))/2 \end{matrix}\right)
\]
depending on wether \( \deg(q_n)<n \) or \( \deg(q_n) =n  \). Thus obtained matrix is still a solution of \hyperref[rhy]{\rhy} and therefore must be equal to \( \boldsymbol Y(z) \). That is, \( q_n(z)=q(z) \) as claimed. The form of the second row can be deduced similarly.

{\bf Step 2.} Let \( J_0 \) be a smooth Jordan curve encircling \( \Delta_0 \) such that \( \rho(z) \) is non-vanishing and analytic on the closure of \( \Omega_0 \), the interior domain of \( J_0 \) with \( \Delta_0 \) removed. We orient \( J_0 \) clockwise when \( \infty\in D_\Delta^0 \) and counter-clockwise otherwise. Further, let \( J_l^+\subset D_\Delta^0 \) and \( J_l^-\subset D_\Delta^\infty \), \( l>0 \), be Jordan curves such that \( \Delta_l \) contains one of them in its interior and another one in its exterior, and \( \rho(z) \) is non-vanishing and analytic on the closures of \( \Omega_l^+ \) and \( \Omega_l^- \), the annular domains with boundaries \( J_l^+\cup\Delta_l \) and \( J_l^-\cup\Delta_l \), respectively. We orient \( J_l^\pm \) in the direction of \( \Delta_l \). Set \( J:=J_0\bigcup\cup_{l>0}(J_l^+\cup J_l^-) \). Assume also that the closures of \( \Omega_0 \) and \( \Omega_l^+\cup\Omega_l^- \), \( l>0 \), are pairwise disjoint.

Recall the definition of \( \varsigma \) in Proposition~\ref{prop:1}. If \( \boldsymbol Y(z) \) is a solution of \hyperref[rhy]{\rhy}, let
\begin{equation}
\label{X}
\boldsymbol X(z) := \boldsymbol Y(z)
\left\{
\begin{array}{rl}
\left(\begin{matrix} 1 & 0 \medskip \\ -\varsigma(v_nw_{\Delta_0}/\rho)(z) & 1 \end{matrix}\right), & z\in\Omega_0, \medskip \\
\left(\begin{matrix} 1 & 0 \medskip \\ \mp (v_nw_{\Delta_0}/\rho)(z) & 1 \end{matrix}\right), & z\in\Omega_i^\pm, \medskip \\
\boldsymbol I, & \text{otherwise}.
\end{array}
\right.
\end{equation}
Then the matrix function \( \boldsymbol X(z) \) solves the following Riemann-Hilbert problem (\rhx):
\begin{myitemize}
\label{rhx}
\item[(a)] \( \boldsymbol X(z) \) is analytic in \( \overline\C\setminus(\Delta\cup J) \) and \( \displaystyle \lim_{z\to\infty} \boldsymbol X(z)z^{-n\sigma_3} = \boldsymbol I \);
\item[(b)] \( \boldsymbol X(z) \) has continuous traces on \( (\Delta\cup J)\setminus\{\pm1\} \) that satisfy 
\[
\boldsymbol X_+(s) = \boldsymbol X_-(s)
\left\{
\begin{array}{rl}
\left(\begin{matrix} 0 & (\rho/(v_nw))(s) \medskip \\ -(v_nw/\rho)(s) & 1 \end{matrix}\right), & s\in\Delta, \medskip \\
\left(\begin{matrix} 1 & 0 \medskip \\ (v_nw_{\Delta_0}/\rho)(s) & 1 \end{matrix}\right), & s\in J;
\end{array}
\right.
\]
\item[(c)] it holds that \( \boldsymbol X(z) = \mathcal O\left(\begin{matrix} 1 & |z-e|^{-1/2} \smallskip \\ 1 & |z-e|^{-1/2}\end{matrix}\right) \) as \( D_\Delta\ni z\to e\in\{\pm1\} \).
\end{myitemize}
In fact, it is quite easy to see that \hyperref[rhx]{\rhx} is solvable if and only if \hyperref[rhy]{\rhy} is solvable and the solutions are connected by \eqref{X} (all matrices in \eqref{X} have determinants identically equal to \( 1 \) and therefore are invertible).

{\bf Step 3.} Let \( \Psi_n(z) \) be given by \eqref{Psin} and \( \varphi(z) = z - \varsigma w_{\Delta_0}(z) \). Define
\begin{equation}
\label{M}
\boldsymbol M(z) := 
\left\{
\begin{array}{rl}
\left(\begin{matrix} \varsigma\Psi_n(z) & \varsigma/(w_{\Delta_0}\Psi_n)(z) \medskip \\ \varsigma(\varphi\Psi_n)(z) & \varsigma/(\varphi w_{\Delta_0}\Psi_n)(z) \end{matrix}\right), & z\in D_\Delta^0, \medskip \\
\left(\begin{matrix} \varsigma\Psi_n(z) & -\varsigma/(w_{\Delta_0}\Psi_n)(z) \medskip \\ \varsigma(\Psi_n/\varphi)(z) & -\varsigma(\varphi/(w_{\Delta_0}\Psi_n))(z) \end{matrix}\right), & z\in D_\Delta^\infty.
\end{array}
\right.
\end{equation}
Using \eqref{P-pr} and the fact that \( \varphi_+(s)\varphi_-(s)\equiv 1 \) for \( s\in\Delta_0 \), one can readily check that \( \boldsymbol M(z) \) solves the following Riemann-Hilbert problem (\rhm):
\begin{myitemize}
\label{rhm}
\item[(a)] \( \boldsymbol M(z) \) is analytic in \( \overline\C\setminus\Delta \) and \( \displaystyle \boldsymbol C^{-1} := \lim_{z\to\infty} \boldsymbol M(z)z^{-n\sigma_3} \) is a diagonal matrix with non-zero entries;
\item[(b)] \( \boldsymbol M(z) \) has continuous traces on \( \Delta\setminus\{\pm1\} \) that satisfy 
\[
\boldsymbol M_+(s) = \boldsymbol M_-(s) \left(\begin{matrix} 0 & (\rho/(v_nw))(s) \medskip \\ -(v_nw/\rho)(s) & 0 \end{matrix}\right), \quad s\in\Delta\setminus\{\pm1\};
\]
\item[(c)] it holds that \( \boldsymbol M(z) =\mathcal O\left(\begin{matrix} 1 & |z-e|^{-1/2} \medskip \\ 1 & |z-e|^{-1/2}\end{matrix}\right) \) as \( D_\Delta\ni z\to e\in\{\pm1\} \).
\end{myitemize}
Observe also that
\begin{equation}
\label{det}
\det(\boldsymbol M(z)) = \frac{1}{(w_{\Delta_0}\varphi)(z)} - \frac{\varphi(z)}{w_{\Delta_0}(z)} \equiv 2\varsigma, \quad z\in \C.
\end{equation}

{\bf Step 4.} Consider the following Riemann-Hilbert problem (\rhz):
\begin{myitemize}
\label{rhz}
\item[(a)] \( \boldsymbol Z(z) \) is analytic in \( \overline\C\setminus J \) and \( \boldsymbol Z(\infty) = \boldsymbol I\);
\item[(b)] \( \boldsymbol Z(z) \) has continuous traces on \( J \) that satisfy 
\[
\boldsymbol Z_+(s) = \boldsymbol Z_-(s)\boldsymbol M(s) \left(\begin{matrix} 1 & 0 \medskip \\ (v_nw_{\Delta_0}/\rho)(s) & 1 \end{matrix}\right)\boldsymbol M^{-1}(s), \quad s\in J.
\]
\end{myitemize}
Denote by \( \boldsymbol J(s) \) the jump matrix in \hyperref[rhz]{\rhz}(b), i.e., \( \boldsymbol Z_+(s) = \boldsymbol Z_-(s)\boldsymbol J(s) \), \( s\in J \). It follows from \eqref{Psin}, \eqref{M}, \eqref{det}, and a straightforward computation that
\[
\boldsymbol J(s) = \boldsymbol I + b_{2n}(s)\frac{\varsigma}{2(w_{\Delta_0}\rho S_\rho^2)(s)}\left(\begin{matrix} \varphi^{-1}(s) & -1 \\ \varphi^{-2}(s) & -\varphi^{-1}(s) \end{matrix}\right) = \boldsymbol I + \boldsymbol o(1)
\]
uniformly on \( J \cap D_\Delta^0 \) as well as
\[
\boldsymbol J(s) = \boldsymbol I + \frac1{b_{2n}(s)}\frac{\varsigma S_\rho^2(s)}{2(w_{\Delta_0}\rho)(s)}\left(\begin{matrix} \varphi(s) & -1 \\ \varphi^2(s) & -\varphi(s) \end{matrix}\right) = \boldsymbol I + \boldsymbol o(1)
\] 
uniformly on \( J \cap D_\Delta^\infty \), where we used Assumption~\ref{as:1}(ii) and the very definition of \( D_\Delta^0 \) and \( D_\Delta^\infty \). The above relations and \cite[Corollary 7.108]{Deift} yield that \hyperref[rhz]{\rhz} is solvable for all \( n \) large enough and the solution satisfies \( \boldsymbol Z(z) = \boldsymbol I + \boldsymbol o(1) \) uniformly in \( \overline \C \).

It now can be quite readily checked that the solution of \hyperref[rhx]{\rhx} is given by \( \boldsymbol X(z) = \boldsymbol C\boldsymbol Z(z)\boldsymbol M(z) \). Observe that \( (1,1) \)-entry of \( \boldsymbol C \) is \( \varsigma\gamma_n \), see \eqref{P-pr} and \eqref{M}. Given a closed set \( K\subset D_\Delta \), the arcs comprising \( J \) can be chosen so that \( K \) does not intersect the closures of \( \Omega_0 \) and \( \Omega_l^+\cup\Omega_l^- \), \( l>0 \). Then it follows from \eqref{X} that \( \boldsymbol X(z) = \boldsymbol Y(z) \) for \( z\in K \). Hence,
\[
q_n(z) = \gamma_n \Psi_n(z)
\left\{
\begin{array}{rl}
1+o(1) + o(1)\varphi(z), & z\in D_\Delta^0, \medskip \\
1+o(1) + o(1)\varphi^{-1}(z), & z\in D_\Delta^\infty,
\end{array}
\right.
\]
uniformly on \( K \). Since \( \varphi^\varsigma(\infty) =0 \), relations \eqref{asymp1} follow. Similarly, \eqref{asymp2} holds since
\[
R_n(z) = \frac{\gamma_n}{w_{\Delta_0}(z)\Psi_n(z)}
\left\{
\begin{array}{rl}
1+o(1) + o(1)\varphi^{-1}(z), & z\in D_\Delta^0, \medskip \\
-\big(1+o(1) + o(1)\varphi(z)\big), & z\in D_\Delta^\infty,
\end{array}
\right.
\]
uniformly on \( K \), where one needs to observe that the error terms \( o(1) \) must vanish at infinity. Asymptotic formula \eqref{asymp3} now follows from \eqref{SKind}, \eqref{Psin}, \eqref{asymp1}, and \eqref{asymp2}.
\end{proof}

\end{document}